\documentclass[11pt]{amsart}

\usepackage[centertags]{amsmath}
\usepackage{amsfonts}
\usepackage{amssymb}
\usepackage{amsthm}
\usepackage[english]{babel}
\usepackage[active]{srcltx}

\usepackage[colorlinks]{hyperref}

\newcommand{\Real}{\mathbb{R}}
\newcommand{\C}{\mathbb{C}}

\newcommand{\co}{\operatorname{co}}

\newcommand{\Lip}{\operatorname{Lip}}
\newcommand{\Spec}{\operatorname{Spec}}
\newcommand{\adj}{\operatorname{adj\,}}
\newcommand{\Isom}{\operatorname{Isom}}

\newcommand{\Ima}{\operatorname{Im}}

\newcommand{\lctan}{\mathcal{L}(T_xM, T_{f(x)}N)}
\newcommand{\derinvpsi}{d\psi^{-1}(\psi(f(x)))}
\newcommand{\derinvphi}{d\varphi^{-1}(\varphi(x))}
\newcommand{\derphi}{d\varphi(x)}
\newcommand{\derpsi}{d\psi(f(x))}

\newcommand{\normJaco}{\|\Jaco\|_{x,\varphi,\psi}}
\newcommand{\conormJaco}{/\hspace{-0.1cm}/ \Jaco /\hspace{-0.1cm}/_{x,\varphi,\psi}}
\newcommand{\conormjax}{/\hspace{-0.1cm}/\hspace{-0.1cm}/ \partial f(x) /\hspace{-0.1cm}/\hspace{-0.1cm}/_{x}}

\newcommand{\conorA}{/\hspace{-0.1cm}/ A /\hspace{-0.1cm}/_{x, \varphi,\psi}}

\newcommand{\Jaco}{\partial (\psi \circ f \circ \varphi^{-1})(\varphi(x))}

\newtheorem{thm}{Theorem}[section]
\newtheorem{cor}[thm]{Corollary}
\newtheorem{lem}[thm]{Lemma}
\newtheorem{prop}[thm]{Proposition}
\newtheorem{properties}[thm]{Properties}

\newtheorem{rem}[thm]{Remark}

\newtheorem{defn}[thm]{Definition}
\newtheorem{notation}[thm]{Notation}

\numberwithin{equation}{section}

\begin{document}

\title{Global inversion of nonsmooth mappings on Finsler manifolds}

\author{Jes\'us A. Jaramillo, \'Oscar  Madiedo and Luis S\'anchez-Gonz\'alez}

\address{Departamento de An{\'a}lisis Matem{\'a}tico\\ Facultad de
Matem{\'a}ticas\\ Universidad Complutense de Madrid\\ 28040 Madrid, Spain}

\thanks{Research supported in part by MICINN Project MTM2009-07848 (Spain). J. A. Jaramillo is also supported by Fundaci\'{o}n Cajamadrid Grants 2011-1012. O. Madiedo is also supported by grant MICINN {BES2010-031192}.  L. S\'anchez-Gonz\'alez has also been supported by grant MEC AP2007-00868.}

\email{jaramil@mat.ucm.es, oscar.reynaldo@mat.ucm.es, lfsanche@mat.ucm.es}


\date{January 18,  2012}

\maketitle

\begin{abstract}
We consider the problem of finding sufficient conditions for a locally Lipschitz mapping  between Finsler manifolds to be a global homeomorphism. For this purpose, we develop the notion of Clarke generalized differential
in this context and, using this, we obtain a version of the Hadamard integral condition for invertibility.
As consequences, we deduce some global inversion and global injectivity  results for Lipschitz mappings on
$\mathbb{R}^n$ in terms of spectral conditions of its Clarke generalized differential.
\end{abstract}

\section{Introduction}

Let $f:{\mathbb R^n} \to {\mathbb R^n}$ be a $C^1$-smooth  mapping with everywhere nonzero Jacobian, so   $f$ is a local diffeomorphism. A natural problem is to study under which conditions we can obtain that $f$ is, in fact, a global diffeomorphism or, equivalently, when  $f$  is globally invertible. This was first considered by Hadamard \cite{Hadamard} who obtained a sufficient condition in terms of the growth of the norm $\Vert [df(x)]^{-1}\Vert$, by means of his celebrated  {\it integral condition}. Namely, $f$ is a global diffeomorphism provided
$$
\int _0^\infty \inf_{\vert\vert x \vert\vert = t}\left\|[df(x)]^{-1}\right\|^{-1} \, dt=\infty.
$$
A wide number of extensions and variants of this result have been obtained in different contexts. In this way, an extension of Hadamard integral condition to the case of local diffeomorphisms between Banach spaces was given by Plastock in \cite{Plastock}. In a nonsmooth setting,  John \cite{John} obtained a variant of  Hadamard integral condition for a local homeomorphism $f$ between Banach spaces, using  the lower scalar Dini derivative of $f$. In the finite-dimensional case,  Pourciau studied in \cite{Pourciau1988}  the global inversion of locally Lipschitz mappings in $\mathbb R^n$ using Clarke generalized jacobian, and he also obtained a variant of Hadamard integral condition in this context. The global invertibility of local diffeomorphisms between Banach-Finsler manifolds has been studied by  Rabier in \cite{Rabier}. For recent global inversion results in a metric space setting, including a version of Hadamard integral condition in terms of an analog of lower scalar derivative, we refer to \cite{GutuJaramillo} and \cite{GaGuJa}.
\noindent On the other hand, and with a different approach, some global invertibility and global injectivity theorems have been recently obtained by  Fernandes,  Guti\'{e}rrez and  Rabanal \cite{FeGuRa} and by  Biasi, Guti\'{e}rrez and dos Santos \cite{BiGuSan}. In particular, they provide conditions of spectral nature on a Lipschitz local diffeomorphism $f$ in $\mathbb R^n$  in order to obtain that $f$ is globally injective or globally invertible.
\smallskip

Our main purpose in this paper is to obtain sufficient conditions of a locally Lipschitz mapping between Finsler manifolds to be a global homeomorphism. Recall that a Finsler manifold consists in a Finsler structure defined on a smooth manifold $M$, i.e., for each point $x\in M$  it assigns  a norm $\Vert \cdot \Vert_x$ in the tangent space $T_xM$ to $M$ at $x$, which depends continuously on $x$. Note that, in particular, every Riemannian manifold is a Finsler manifold. Our global inversion results are given  in terms of an analog of Clarke generalized jacobian, extending the results of  Pourciau \cite{Pourciau1988} to this setting.  Furthermore, we  also give an extension of the results of \cite{FeGuRa} and \cite{BiGuSan} to the context of nonsmooth Lipschitz mappings in $\mathbb R^n$.

\noindent The paper is organized as follows. In Section 2 we introduce the notion of Clarke generalized differential in the general setting of  smooth manifolds, and  some fairly basic properties are obtained. In Section 3, we focus on the case of Finsler manifolds. The norm of Clarke generalized differential is introduced, and we explore its connections with the upper scalar Dini derivative and with the local Lipschitz constant. Section 4 is devoted to global inversion theorems. In the case of   locally Lipschitz mappings between Finsler manifolds, we obtain a version of Hadamard integral condition in terms of the co-norm of Clarke generalized differential. In the case of  Lipschitz mappings in $\mathbb R^n$, we obtain a variant of Hadamard integral condition involving the eigenvalues of Clarke generalized differential.

\bigskip
\section{Generalized differential on smooth manifolds}

Firstly, let us recall the definition of Clarke generalized differential of a mapping $f:\Real^m \to \Real^n$ which is Lipschitz on a neighborhood of $x\in \Real^m$. In this case, by the classical Rademacher theorem, $f$ is almost everywhere differentiable in a neighborhood of $x$, and according to \cite{Clarkebook}  the \emph{Clarke generalized differential} of $f$ at $x$ is defined as
\begin{equation*}
\partial f(x) := \co \, \left\{\lim_{i\to\infty} df(x_i): x_i\to x \ \text{and}\ df(x_i) \text{ exists}\right\}.
\end{equation*}
Here $df(x)$ denotes the differential of $f$ at $x$, which we identify with the corresponding jacobian matrix, and {\it co} means, as usual, the convex hull. It can be proved that $\partial f(x)$ is then a nonempty, convex and compact set of  linear mappings from $\Real^m $ to $ \Real^n$ (see \cite[Proposition 2.6.2]{Clarkebook}). We refer to \cite{Clarkebook, Clarke1975, Pourciau1977} for further information and basic properties of generalized differentials.

\medskip

 We extend this definition  to the context of  smooth manifolds. Our notation here is standard. We refer to \cite{Lang, Deim}  for basic definitions and notations on manifolds. If $E$ and $F$ are vector spaces, we denote by $\mathcal{L}(E, F)$ the space of all linear mappings from $E$ to $F$.

\begin{defn}\label{defLocalLips}
Let $M$ and $N$ be $C^1$-smooth manifolds of dimension $m$ and $n$, respectively. A mapping $f:M\to N$ is said to be Lipschitz on a neighborhood of $x\in M$ if there are  charts $(U,\varphi)$ of $M$ at $x$ and $(V,\psi)$  of $N$ at $f(x)$,  with $f(U)\subset V$  such that $\psi \circ f \circ \varphi^{-1} : \varphi(U)\subset \Real^m \to \Real^n$ is Lipschitz on a neighborhood of $\varphi(x)$. A mapping $f:M\to N$ is said to be locally Lipschitz if  it is Lipschitz on a neighborhood of  $x$ for every $x\in M$.
\end{defn}

Notice that every $C^1$-smooth mapping between euclidean spaces is locally Lipschitz. Thus the above definitions do not depend on charts, and we also have that every $C^1$-smooth mapping between $C^1$-smooth manifolds is locally Lipschitz.

\begin{defn}
Let $M$ and $N$ be $C^1$-smooth manifolds, and let $f:M\to N$ be a map which is Lipschitz on a neighborhood of  $x\in M$. Then, the \textbf{Clarke generalized differential} of $f$ at $x$ is defined as
\begin{align*}
\partial f(x) :&= \{ d\psi^{-1}(\psi(f(x)))\circ A \circ d\varphi(x): A\in \partial (\psi \circ f \circ \varphi^{-1}) (\varphi(x))\} \\
& = d\psi^{-1}(\psi(f(x)))\circ \partial (\psi\circ f \circ \varphi^{-1})(\varphi(x))\circ d\varphi(x),
\end{align*}
where $(U,\varphi)$ is a chart of $M$ at $x$ and $(V,\psi)$ is a chart of $N$ at $f(x)$. Notice that $\partial f(x)\subset \lctan$.
\end{defn}

\begin{prop}
The above definition does not depend on charts.
\end{prop}
\begin{proof}
Let us take $(U_1,\varphi_1)$, $(U_2,\varphi_2)$ different charts of $M$ at $x$ and $(V_1,\psi_1)$, $(V_2,\psi_2)$ different charts of $N$ at $f(x)$. Let us denote by
\begin{align*}
\partial_1 f(x): &= \{ d\psi_1^{-1}(\psi_1(f(x)))\circ A \circ d\varphi_1(x): A\in \partial (\psi_1 \circ f \circ \varphi_1^{-1}) (\varphi_1(x))\}, \text{ and}\\
\partial_2 f(x):&= \{ d\psi_2^{-1}(\psi_2(f(x)))\circ A \circ d\varphi_2(x): A\in \partial (\psi_2 \circ f \circ \varphi_2^{-1}) (\varphi_2(x))\}.
\end{align*}
We have to show that $\partial_1 f(x) = \partial_2 f(x)$. Let us note that $\partial_1 f(x)$ and $\partial_2 f(x)$ are convex sets. Thus, it is sufficient to prove that
\begin{align*}
B\in \partial_2 f(x) \ \text{ for all } \ & B=d\psi_1^{-1}(\psi_1(f(x)))\circ A \circ d\varphi_1(x)  \in \partial_1 f(x)\ \text{ with}\\
&A=\lim\limits_{i\to \infty} d(\psi_1 \circ f \circ \varphi_1^{-1})(\varphi_1(x_i)) \in \partial (\psi_1 \circ f \circ \varphi_1^{-1}) (\varphi_1(x)),
\end{align*}
where $x_i \to x$ and $d(\psi_1 \circ f \circ \varphi_1^{-1})(\varphi_1(x_i))$ exists.

Let us take $B$ and $A$ with the above hypothesis. Since $\psi_2\circ \psi_1^{-1}$ and $\varphi_1Ê\circ \varphi_2^{-1}$  are $C^1$-smooth in a neighborhood of $\psi_1(f(x))$ and   in a neighborhood of $\varphi_2(x)$, respectively, then
\begin{equation*}
(\psi_2\circ \psi_1^{-1})  \circ      (\psi_1 \circ f \circ \varphi_1^{-1}) \circ     (\varphi_1Ê\circ \varphi_2^{-1}) = \psi_2 \circ f \circ \varphi_2^{-1}
    \end{equation*}
is differentiable at $\varphi_2(x_i)$ for every $i$ large enough. Moreover,
\begin{align*}
\lim_{i\to \infty}  d(\psi_2 \circ f \circ \varphi_2^{-1}&)(\varphi_2(x_i)) \\
& = \lim_{i\to \infty} d \left((\psi_2\circ \psi_1^{-1})  \circ      (\psi_1 \circ f \circ \varphi_1^{-1})     \circ     (\varphi_1Ê\circ \varphi_2^{-1})\right) (\varphi_2(x_i))     \\
& = d(\psi_2\circ \psi_1^{-1}) (\psi_1(f(x))) \circ A \circ d(\varphi_1Ê\circ \varphi_2^{-1})(\varphi_2(x)) Ê\\
& =  d\psi_2 (f(x)) \circ B \circ d\varphi_2^{-1}(\varphi_2(x)) = C \in \partial (\psi_2 \circ f \circ \varphi_2^{-1})(\varphi_2(x)).
\end{align*}
Thus,
$B=d\psi_2^{-1}(\psi_2(f(x)))\circ C \circ d\varphi_2(x)$  with  $C\in \partial (\psi_2 \circ f \circ \varphi_2^{-1}) (\varphi_2(x))$, and $B\in \partial_2 f(x)$.
\end{proof}

\

\begin{rem}\label{defbundle}
{\rm Let $M$ and  $N$ be $C^1$-smooth manifolds of dimension $m$ and $n$, respectively, and $f:M\to N$ a continuous mapping. We consider the vector bundle $\mathcal{L}(TM,f^*TN)$ defined as the following disjoint union:
\begin{equation*}
\mathcal{L}(TM,f^*TN) :=\bigsqcup_{x\in M}\lctan =\{(x,B):x\in M \text{ and } B \in \lctan \},
\end{equation*}
where the bundle projection $\pi:\mathcal{L}(TM,f^*TN)\to M$ is given by $\pi(x,B)=x$ (we refer for instance to \cite{Milnor} for details).

Recall that the topology of $\mathcal{L}(TM,f^*TN)$ is generated by the basic open sets ${\mathcal O}(U, V, W)$ described below, where $(U, \varphi)$ is a chart of $M$, $(V, \psi)$ is a chart of $N$ with $f(U)\subset V$ and $W$ is an open subset of $\mathcal{L}(\Real^m,\Real^n)$. Here we denote
\begin{equation*}
{\mathcal O}(U, V, W)=\bigsqcup_{x\in U} W_{x, \varphi, \psi}=\{(x,B):x\in U \text{ and } B \in W_{x, \varphi, \psi}\}
\end{equation*}
and $W_{x, \varphi, \psi}=\{\derinvpsi\circ A \circ \derphi: A\in W \}=\derinvpsi\circ W \circ \derphi$  for every $x\in U$.

Thus, if $(U, \varphi)$ is a chart of $M$ with $x \in U$ and $(V, \psi)$ is a chart of $N$ with $f(U)\subset V$, we obtain a local trivialization of $\mathcal{L}(TM,f^*TN)$ around $x$ by means of the homeomorphism
$$
h: \pi^{-1} (U)= {\mathcal O}(U, V, \mathcal{L}(\Real^m,\Real^n)) \rightarrow U \times \mathcal{L}(\Real^m,\Real^n)
$$
given by
$$
h (x, B)= \left(x, \, d\psi(f(x))\circ B \circ d\varphi^{-1}(\varphi(x))\right).
$$
}
\end{rem}

Now, we list some basic properties of Clarke generalized differential on manifolds.

\begin{properties}\label{basic:properties}
Let $M$ and $N$ be $C^1$-smooth manifolds. Let $f:M\to N$ be a mapping which is Lipschitz on a neighborhood of $x\in M$. Then:
\begin{enumerate}
\item $\partial f(x)$ is a nonempty, convex and compact subset of $\lctan$.
\item If the sequence $(x_i, B_i)$ converges to $(x,B)$ in $\mathcal{L}(TM, f^*TN)$ and  $B_i\in \partial f(x_i)$ for each $i$, then  $B\in \partial f(x)$.
\item If $f$ is a $C^1$-smooth mapping then $\partial f(x)=\{df(x)\}$.
\end{enumerate}
\end{properties}
\begin{proof}
It is straightforward to check these properties. For reader's convenience, we shall show  property $(2)$.

Let us take $(U,\varphi)$ and $(V,\psi)$ charts of $M$ at $x$ and $N$ at $f(x)$, respectively, with $f(U)\subset V$, and consider the local trivialization of $\mathcal{L}(TM,f^*TN)$ around $x$ given by
$$
h: \pi^{-1} (U) \rightarrow U \times \mathcal{L}(\Real^m,\Real^n)
$$
where
$$
h (x, B)= \left(x, \, d\psi(f(x))\circ B \circ d\varphi^{-1}(\varphi(x))\right).
$$
For $i$ large enough we have that  $x_i \in U$ and $f(x_i) \in V$. In this case, let us denote $A_i= d\psi(f(x_i)) \circ B_i \circ d\varphi^{-1}(\varphi(x_i)) \in \partial(\psi \circ f \circ \varphi^{-1})(\varphi(x_i))\subset \mathcal{L}(\Real^m,\Real^n)$. By continuity, we have that the sequence $(A_i)$  converges  to
$$
 A:= d\psi(f(x))\circ B \circ d\varphi^{-1}(\varphi(x))  \in \mathcal{L}(\Real^m,\Real^n).
$$
Now by \cite[Proposition 2.6.2(b)]{Clarkebook} we obtain that $A\in \partial(\psi \circ f\circ \varphi^{-1})(\varphi(x))$,
and therefore $B= d\psi^{-1}(\psi(f(x)))\circ A \circ d\varphi(x) \in \partial f(x)$.
\end{proof}

Next we are going to see some continuity property of the Clarke generalized differential.

\begin{defn}
Let $X$ and $Y$ be topological spaces. A   multivalued mapping $\Phi:X\to 2^Y$ is  upper semicontinuous (usc) at $x\in X$ if for any open set ${\mathcal O}$ of $Y$ such that $\Phi(x) \subset {\mathcal O}$, there is an open neighborhood $U$ of $x$ in $X$  such that $\Phi(y)\subset {\mathcal O}$ for every $y\in U$.  A multivalued mapping $\Phi:X\to 2^Y$ is said to be upper semicontinuous (usc)  if it is usc at $x$ for every $x\in X$.
\end{defn}

\begin{prop}
Let $M$ and $N$ be $C^1$-smooth manifolds and let $f:M\to N$ be a locally Lipschitz mapping. Then, the Clarke generalized differential of $f$ is an usc mapping, in the sense that the multivalued mapping
\begin{equation*}
\widetilde{\partial} f : M \to 2^{\mathcal{L}(TM, f^*TN)}
\end{equation*}
given by $\widetilde{\partial}f(x)=(x,\partial f(x))$ is  upper semicontinuous.
\end{prop}
\begin{proof}
Let  $x\in M$ and consider an open subset ${\mathcal O}$ of $\mathcal{L}(TM, f^*TN)$ with $\widetilde{\partial} f (x) \subset  {\mathcal O}$. According to Remark \ref{defbundle}, for each $B \in \partial f(x)$ we can choose a basic open set of the form ${\mathcal O}(U, V, W)$ such that $(x, B)\in {\mathcal O}(U, V, W) \subset {\mathcal O}$. Since $\widetilde{\partial} f (x)$ is compact, it can be covered by a finite family of such basic open sets. Therefore we can assume that ${\mathcal O}={\mathcal O}(U, V, W)$ is a basic open set, where $(U, \varphi)$ is a chart of $M$ with $x\in U$, $(V, \psi)$ is a chart of $N$ with $f(U)\subset V$ and $W$ is an open subset of $\mathcal{L}(\Real^m, \Real^n)$.
Thus  $\partial f(x) \subset W_{x, \varphi, \psi}=\derinvpsi\circ W \circ \derphi$ and
\begin{equation*}
\partial(\psi \circ f \circ \varphi^{-1}) (\varphi(x)) \subset \derpsi \circ W_{x,\varphi,\psi}\circ \derinvphi =W.
\end{equation*}
Since
$$
\partial(\psi \circ f \circ \varphi^{-1}): \varphi (U) \rightarrow 2^{\mathcal{L}(\Real^m, \Real^n)}
$$
is upper semicontinuous at $\varphi(x)$ (see \cite[Proposition 2.6.2(c)]{Clarkebook}), there is an open neighborhood $\widetilde{U}$ of $x$ such that $\partial(\psi \circ f \circ \varphi^{-1}) (\varphi(y))\subset W$ for every $y\in \widetilde{U}$. Then,
by the definition of the Clarke generalized differential, we have that
\begin{equation*}
\partial f(y) \subset  d\psi^{-1}(\psi(f(y)))\circ W\circ d\varphi(y)=W_{y,\varphi,\psi},
\end{equation*}
for every $y\in \widetilde{U}$. Thus, $\widetilde{\partial} f(y)\subset {\mathcal O}(U, V, W)$  for every $y\in\widetilde{U}$. This shows that $f$ is usc at $x$.
\end{proof}

\begin{rem}
Notice that by Property \ref{basic:properties}(2) the multivalued mapping $\widetilde{\partial} f$ is closed at $x$, i.e., if the sequence  $(x,B_i)$ converges to $(x, B)$ in $\mathcal{L}(TM, f^*TN)$, and $(x,B_i)\in  \widetilde{\partial}f(x)$  for every $i$,  then  $(x, B) \in \widetilde{\partial}f(x)$.
\end{rem}


In the following proposition we obtain a chain rule for the Clarke generalized differential.

\begin{prop}
Let $M$ and $N$ be $C^1$-smooth manifolds. Let $h:M\to N$ be a Lipschitz mapping on a neighborhood of $x$ and $g: N \to \Real$ a Lipschitz function on a neighborhood of $h(x)$.  Then $f=g\circ h:M\to \Real$ is Lipschitz on a neighborhood of $x$ and
$$
\partial f(x) \subset \co\{\partial g(h(x)) \circ \partial h(x)\}.
$$
If, in addition, $g$ is $C^1$-smooth, then
$$
\partial f(x) = d g(h(x)) \circ \partial h(x).
$$
\end{prop}
\begin{proof}
We start with the first inclusion. Let us take $(U,\varphi)$ and $(V,\psi)$  charts of $M$ at $x$ and $N$ at $h(x)$,
respectively, with $h(U) \subset V$. By convexity, it is sufficient to prove that
\begin{equation*}
B\in \partial g(h(x)) \circ \partial h(x) \ \text{ for all } \ B=A  \circ d\varphi(x) \ \text{ with }\  A \in \partial (f \circ \varphi^{-1})(\varphi(x)).
\end{equation*}
Let us take $B$ and $A$ with the above hypothesis. Since
\begin{align*}
\partial (f \circ \varphi^{-1})(\varphi(x)) &= \partial (g \circ h \circ \varphi^{-1})(\varphi(x)) \\
&= \partial ((g \circ \psi^{-1})\circ (\psi \circ h \circ \varphi^{-1}))(\varphi(x)),
\end{align*}
and the Clarke generalized differential in $\Real^n$ satisfies the corresponding chain rule
\cite[Theorem 2.6.6]{Clarkebook}, we obtain that
\begin{align*}
 A \in \partial (f \circ \varphi^{-1})&(\varphi(x))=\partial ((g \circ \psi^{-1})\circ (\psi \circ h \circ \varphi^{-1}))(\varphi(x))\\
& \subset \co\{\partial (g \circ \psi^{-1})(\psi(h(x)))\circ 	\partial(\psi \circ h \circ \varphi^{-1})(\varphi(x))\}.
\end{align*}

Without loss of generality we may suppose that $A=A_1\circ A_2$ with  $A_{2} \in \partial(\psi \circ h \circ \varphi^{-1})(\varphi(x))$ and $A_{1} \in \partial(g \circ \psi^{-1})(\psi(h(x)))$. Then
\begin{align*}
B & = A\circ d\varphi(x)  \\
& = \left( A_1\circ  d\psi(h(x)) \right) \circ \left(d\psi^{-1}(\psi(h(x))) \circ A_{2} \circ \derphi \right) \in \partial g(h(x))\circ \partial h(x),
\end{align*}
since $A_{2} \in \partial(\psi \circ h \circ \varphi^{-1})(\varphi(x))$ and $A_{1} \in \partial(g \circ \psi^{-1})(\psi(h(x)))$.

The last assertion follows also from \cite[Theorem 2.6.6]{Clarkebook} in an analogous way.
\end{proof}

Now using the same ideas as above and the result of \cite[Corollary, p. 75]{Clarkebook}, the following property can be proved.

\begin{cor}
Let $M,N$ and $S$ be $C^1$-smooth manifolds. Let $h:M\to N$ be a Lipschitz mapping on a neighborhood of $x$ and $g: N \to S$  a Lipschitz mapping on a neighborhood of $h(x)$. Then $f=g\circ h:M\to S$ is Lipschitz on a neighborhood of $x$ and
$$
\partial f (x) v \subset \co\{\partial g(h(x)) \circ \partial h(x) v\},
$$
for every $v\in T_xM$. If, in addition, $g$ is $C^1$-smooth, then
$$
\partial f(x) = d g(h(x)) \circ \partial h(x) v.
$$
\end{cor}

\section{Generalized differential on Finsler manifolds}

In this section we are going to show some special features of Clarke generalized  differential for locally Lipschitz mappings between Finsler manifolds. Recall that the tangent bundle of a smooth manifold $M$ is $TM = \{(x,v) : x \in M \text{ and } v \in T_xM\}$.

\begin{defn} \label{defFinsler}
A \textbf{$C^1$ Finsler manifold} is a pair $(M,||\cdot||_M)$, where $M$ is a $C^1$-smooth manifold, and $||\cdot||_M: TM\to [0,\infty)$ is a continuous mapping defined on the tangent bundle $TM$ of $M$ which satisfies that
for every $x\in M$, the restriction  $||\cdot||_x:={||\cdot||_M}_{\mid_{T_xM}}:T_xM\to [0,\infty)$ is a norm on the
tangent space $T_xM$.

\end{defn}

\begin{rem}\label{Finsler-Palais}
{\rm Notice that the concept of $C^1$ Finsler manifold coincides with the concept of finite-dimensional $C^1$ Finsler manifold in the sense of Palais \cite{Palais}. This can be proved using the local compactness of the manifold. In particular, if $(M, \|\cdot\|_M)$ is a $C^1$ Finsler manifold of dimension $n$, then for every $x\in M$, every $\varepsilon>0$ and every chart $\varphi:U\to \Real^n$ with $x\in U\subset M$, there is an open neighborhood $W$ with $x\in W\subset U$ such that
\begin{equation}
\frac{1}{1+\varepsilon}||d\varphi^{-1}(\varphi(x))(v)||_{x}\le ||d\varphi^{-1}(\varphi(y))(v)||_{y}\le (1+\varepsilon)||d\varphi^{-1}(\varphi(x))(v)||_{x},
\end{equation}
 for every  $y\in W$ and every $v\in \Real^n$. In terms of equivalence of norms, the above inequalities yield the fact that the norms $||d\varphi^{-1}(\varphi(x))(\cdot)||_{x}$
and   $||d\varphi^{-1}(\varphi(y))(\cdot)||_{y}$ are $(1+\varepsilon)$-equivalent.}
\end{rem}

For further information about Finsler manifolds we refer to  \cite{Palais, Deim, Rabier, MarLuis}. Let us note that, in particular, every Riemannian manifold is a $C^1$ Finsler manifold \cite{Palais}.

If $f:M \to N$ is a $C^1$-smooth mapping between Finsler manifolds, the norm of its differential at the point $x\in M$
is defined  as
\begin{equation*}
||df(x)||_x= \sup\{||df(x)(v)||_{f(x)}:v\in T_x M, ||v||_x = 1\}.
\end{equation*}

The Finsler structure allows us to define  the norm and the co-norm for the Clarke generalized differential of a locally Lipschitz mapping. Let us  recall that  Pourciau \cite{Pourciau1988} considered the co-norm of the Clarke generalized differential in order to obtain a formulation of Hadamard integral condition in a nonsmooth setting. According to \cite{Pourciau1988}, the co-norm of a linear mapping $A:\Real^m \to \Real^n$ is defined as
\begin{equation*}
/\hspace{-0.1cm}/A/\hspace{-0.1cm}/ = \inf_{||u||=1}||Au||.
\end{equation*}
Therefore, in an analogous way,  we can now introduce the following definitions.

\begin{defn}
Let $M$ and $N$ be $C^1$ Finsler manifolds, and let $f: M \to N$ be a Lipschitz mapping on a neighborhood of
$x\in M$. The norm and the co-norm of the Clarke generalized differential of $f$ at $x$ are defined, respectively, as
\begin{align*}
|||\partial f(x)|||_x & =  \sup\{|||B|||_x: B \in \partial f(x)\} \quad \text{and} \\
 \conormjax  & =  \inf\{/\hspace{-0.1cm}/\hspace{-0.1cm}/ B /\hspace{-0.1cm}/\hspace{-0.1cm}/_x: B \in \partial f(x)\},
\end{align*}
where for each $B\in \lctan$, the norm and the co-norm of $B$ are defined by
\begin{equation*}
|||B|||_{x} = \sup\limits_{\substack{v \in T_xM \\
                                          ||v||_x =1 }} ||Bv||_{f(x)}\quad \text{and} \quad
                              /\hspace{-0.1cm}/\hspace{-0.1cm}/ B /\hspace{-0.1cm}/\hspace{-0.1cm}/_{x} = \inf\limits_{\substack{v \in T_xM \\
                                          ||v||_x =1 }} ||Bv||_{f(x)}.
\end{equation*}
\end{defn}

The following notation will be helpful along this paper:

\begin{notation}
Let $M$ and $N$ be $C^1$ Finsler manifolds of dimension $m$ and $n$, respectively.  Let $f:M\to N$ be a mapping, let $x \in M$, and consider charts $(U,\varphi)$ of  $M$ at $x$ and $(V,\psi)$ of $N$ at $f(x)$.  We denote  by
\begin{itemize}
\item $\|\cdot\|_{x, \varphi}: \Real^{m} \to [0,\infty)$  the  norm $\|d\varphi^{-1}(\varphi(x))(\cdot)\|_{x}$ on $\Real^{m}$.
\vspace{0.2cm}

\item $\|\cdot\|_{f(x), \psi}: \Real^{n} \to [0,\infty)$  the  norm $\|d\psi^{-1}(\psi(f(x)))(\cdot)\|_{f(x)}$ on $\Real^{n}$.
\vspace{0.2cm}
\item $\|\cdot\|_{x,\varphi, \psi}: \mathcal{L}(\Real^{m},\Real^{n}) \to [0,\infty)$ the norm $\|A\|_{x,  \varphi, \psi} = \sup\{\|Av\|_{f(x), \psi}: \|v\|_{x,\varphi} = 1\}$ on $\mathcal{L}(\Real^{m},\Real^{n})$.
\vspace{0.2cm}
\item $/\hspace{-0.1cm}/\cdot/\hspace{-0.1cm}/_{x,\varphi,\psi}: \mathcal{L}(\Real^{m},\Real^{n})  \to [0,\infty)$ the co-norm $\conorA = \inf\{||Av||_{f(x),\psi}: ||v||_{x,\varphi} = 1\}$ on $ \mathcal{L}(\Real^{m},\Real^{n})$.
\end{itemize}
\end{notation}

Our next result shows how the norm and the co-norm of the Clarke generalized differential can be obtained via localization through arbitrary charts.

\begin{prop}\label{co-norm:igualdad}
Let $M$ and $N$ be $C^1$ Finsler manifolds and $f: M \to N$ a Lipschitz mapping on a neighborhood of $x$. Then
\begin{align*}
|||\partial f(x)|||_x & = \normJaco \quad \text{and}\\
\conormjax & = \conormJaco
\end{align*}
for every pair of charts $(U,\varphi)$ and $(V,\psi)$ of $M$ at $x$ and $N$ at $f(x)$, respectively.
\end{prop}

\begin{proof}
The proof  is immediate from the definitions and the following lemma:

\begin{lem}\label{norm:and:conorm}
Under the assumptions of Proposition \ref{co-norm:igualdad},
let $B \in \lctan$ which can be written by $B = \derinvpsi \circ A \circ \derphi$ with $A\in \mathcal{L}(\Real^m,\Real^n)$. Then
\begin{itemize}
\item[(a)] $|||B|||_x = \|A\|_{x,\varphi,\psi}$, and
\item[(b)] $/\hspace{-0.1cm}/\hspace{-0.1cm}/ B /\hspace{-0.1cm}/\hspace{-0.1cm}/ _x= \conorA$.
\end{itemize}
\end{lem}

Let us prove the first equality,  the proof of (b) will follow along the same lines.

\begin{align*}
|||B|||_x &= \sup\limits_{{\|v\|}_{x}=1}\|Bv\|_{f(x)}  = \sup\limits_{{\|v\|}_{x}=1}\|(\derinvpsi \circ A \circ \derphi)v\|_{f(x)}  \\
& = \sup\limits_{{\|v\|}_{x}=1}\|A \circ \underbrace{\derphi v}_{w}\|_{f(x),\psi} = \sup\limits_{\|\derinvphi w\|_{x}=1}\|A w\|_{f(x),\psi} \\
&= \sup\limits_{\|w\|_{x,\varphi}=1}\|Aw\|_{f(x),\psi} =  \|A\|_{x,\varphi,\psi}.
\end{align*}
\end{proof}

Every {\it connected} Finsler manifold  $M$ has a natural metric structure, which is defined in the following way. First recall that the {\em length} of a piecewise $C^1$-smooth
path $\gamma:[a,b]\rightarrow M$ is defined as
\begin{equation}\label{lenght}
\ell(\gamma):=\int_{a}^b||\gamma'(t)||_{\gamma(t)}\,dt.
\end{equation}
Now since $M$ is connected, it is connected by piecewise $C^1$-smooth paths, and the associated {\em Finsler distance} $d_M$ on $M$ is defined as
\begin{equation*}
d_M(x,y)=\inf\{\ell(\gamma): \, \gamma \text{ is a piecewise } C^1 \text{-smooth path connecting } x \text{ and } y\}.
\end{equation*}
This Finsler metric is consistent with the manifold topology given in $M$ (see \cite{Palais}). The open ball and the closed ball of center $x\in M$ and radius $r>0$ will be denoted by $B_M(x,r):=\{y\in M:\, d_M(x,y)<r\}$ and $\overline{B}_M(x,r):=\{y\in M:\, d_M(x,y)\le r\}$, respectively.

If $M$ and $N$ are connected Finsler manifolds, we say that a mapping $f:M \rightarrow N$ is {Lipschitz} if  $f:(M, d_M)\rightarrow (N, d_N)$ is Lipschitz for the corresponding Finsler distances. Then, its Lipschitz constant $\Lip(f)$ is defined, as usual, as
$$
\Lip(f)=\sup\left\{\frac{d_N(f(x),f(y))}{d_M(x,y)}: x,y \in M, x\not=y\right\}.
$$

Let us now recall the following mean value inequality for Finsler manifolds from \cite[Proposition 2.3]{MarLuis} (see also \cite{AzFeMe}).

\begin{lem}  \label{meanvalue} Let $M$ and $N$ be connected $C^1$ Finsler manifolds
and $f:M\to N$  a  $C^1$-smooth  mapping. Then,  $f$ is Lipschitz if and only if  $||df||_\infty := \sup\{||df(x)||_x \, : \, x\in M\}<\infty$. Furthermore, $\Lip(f)= ||df||_\infty$.
\end{lem}

We will also need the following result related to the  $(1+\varepsilon)$-bi-Lipschitz local behavior of charts of a $C^1$ Finsler manifold. This follows at once from \cite[Lemma 2.4]{MarLuis}, taking into account Remark \ref{Finsler-Palais}.

\begin{lem} \label{desigualdades:BiLipschitz}
Let us consider a connected $C^1$ Finsler manifold $M$. Then, for every $x\in M$, every chart
$(U,\varphi)$  with $x\in U$ and every $\varepsilon>0$,  there exists an open neighborhood $W\subset U$ of $x$ satisfying
\begin{equation}
\frac{1}{1+\varepsilon}d_M(y, z)\le \|\varphi(y)-\varphi(z)\|_{x,\varphi}\le (1+\varepsilon) d_M(y, z), \quad \text{ for every } y, z\in W.
\end{equation}
\end{lem}

\begin{rem}
{\rm As a consequence of the above lemma, it is easy to see that if $M$ and $N$ are connected $C^1$ Finsler manifolds, a mapping $f:M\to N$  is locally Lipschitz in the sense of Definition \ref{defLocalLips} if, and only if, $f:(M, d_M)\to (N, d_N)$ is locally Lipschitz for the corresponding Finsler metrics.
}
\end{rem}

Let us give the definition of the lower and upper scalar Dini derivatives of a continuous mapping between Finsler manifolds. These quantities were used by  John  \cite{John} for local homeomorphisms between Banach spaces in order to obtain a version of Hadamard integral condition. Later on, this kind of scalar derivatives have been considered in \cite{GutuJaramillo} and \cite{GaGuJa}, in a metric space setting.

\begin{defn}
Let $M$ and $N$ be connected $C^1$ Finsler manifolds, $f:M\to N$  a  continuous mapping and $x\in M$.
The \textbf{lower and upper scalar derivatives} of $f$ at $x$ are defined as
\begin{equation*}
D^{-}_x f =\liminf_{y\to x} \frac{d_N(f(y),f(x))}{d_M(y,x)}, \quad D^{+}_x f =\limsup_{y\to x} \frac{d_N(f(y),f(x))}{d_M(y,x)}
\end{equation*}
where $y\in M$ and $y\neq x$.
\end{defn}

It was proved by  John \cite{John} that, if $E$ and $F$ are Banach spaces and $f:E\to F$ is differentiable at $x\in E$, then $D^+_xf=||df(x)||$ and, if in addition $df(x)$ is invertible, $D^-_xf=||[df(x)]^{-1}||^{-1}$. An analogous statement holds for smooth mappings between Riemannian manifolds: If $f : M \to N$ is a $C^1$-smooth mapping between connected and complete Riemannian manifolds, then $D^+_x f = ||df(x)||_x$ and, if in addition $df(x)\in \Isom(T_xM, T_{f(x)}N)$, then $D^{-}_x f = ||[df(x)]^{-1}||_{f(x)}^{-1}$ (see \cite{GutuJaramillo}). We shall  prove that the same property holds for connected and complete  Finsler manifolds. Let us recall that a Finsler manifold $M$ is said to be \emph{complete} if it is complete with respect to its Finsler metric $d_M$.

\begin{prop}
Let $f:M\to N$ be a $C^1$-smooth mapping between  connected and complete $C^1$ Finsler manifolds. Then,   we have that $D^+_x f = ||df(x)||_x$ for every $x\in M$. If in addition $d f(x)\in \Isom(T_x M, T_{f(x)}N)$, then $D^-_x f =||[df(x)]^{-1}||_{f(x)}^{-1}$.
\end{prop}
\begin{proof}
From \cite[Example 3.3]{GutuJaramillo}, we have that $D^+_x f \le ||df(x)||_x$, and if  $df(x)\in \Isom(T_xM, T_{f(x)}N)$ then $D^{-}_x f \ge ||[df(x)]^{-1}||_{f(x)}^{-1}$. So, it only remains to prove the reverse inequality.

Let us take $\varepsilon>0$. By  Lemma \ref{desigualdades:BiLipschitz}, we can find a chart $(U,\varphi)$ of $M$ with  $x\in U$ and a chart $(V,\psi)$ of $N$ with $f(U) \subset V$ such that $\varphi(x)=0$,  $\psi(f(x))=0$ and the charts  $\varphi$ and $\psi$ are $(1+\varepsilon)$-bi-Lipschitz for the norms $||\cdot||_{x,\varphi}$ on $\mathbb R^m$  and $||\cdot||_{f(x),\psi}$ on $\mathbb R^n$, respectively. Fix $v\in T_xM$ with $||v||_x=1$, and denote $w=d\varphi (x) (v) \in \Real^m$; then $||w||_{x,\varphi}=1$. Choose some  $\delta>0$ such that $tw\in \varphi(U)$ for all $t\in (-\delta,\delta)$, and consider the $C^1$-smooth curve $\gamma:(-\delta,\delta)\to M$ given by $\gamma(t)=\varphi^{-1}(tw)$, for
every $t\in (-\delta,\delta)$. Using that $\psi$ is $(1+\varepsilon)$-Lipschitz for the norm $||\cdot||_{f(x),\psi}$, we obtain that
\begin{align*}
||(\psi\circ f\circ \gamma)'(0)||_{f(x),\psi}= & \left\|\lim_{t\to 0} \frac{\psi\circ f\circ \gamma(t)-\psi\circ f\circ \gamma(0)}{t}\right\|_{f(x),\psi} \\
\le & (1+\varepsilon)\lim_{t\to 0} \frac{d_N(f(\gamma(t)), f(x))}{|t|}.
\end{align*}
Now, since $\varphi^{-1}$ is $(1+\varepsilon)$-Lipschitz for  the norm $||\cdot||_{x,\varphi}$,
we can  assure that
\begin{equation*}
d_M(\gamma(t),x)=d_M(\varphi^{-1}(tw),\varphi^{-1}(0))\le {(1+\varepsilon)}||tw-0||_{x,\varphi}={(1+\varepsilon)}|t|.
\end{equation*}
Hence, since $\gamma'(0)=d\varphi^{-1}(0)(w)=v$,
\begin{align*}
||df(x)(v)||_{f(x)}=&||d(f\circ \gamma)(0)||_{f(x)}=||d(\psi^{-1}\circÊ\psi \circ f\circ \gamma)(0)||_{f(x)}\\
=& ||d\psi^{-1}(0)(\psi\circ f\circ \gamma)'(0)||_{f(x)}=||(\psi\circ f\circ \gamma)'(0)||_{f(x),\psi} \\
\le & (1+\varepsilon)^2 \lim_{t\to 0} \frac{d_N(f(\gamma(t)),f(x))}{d_M(\gamma(t),x)}
\le (1+\varepsilon)^2D_x^+f.
\end{align*}
This holds for  every  $v\in T_xM$ with $||v||_x=1$ and every $\varepsilon>0$, so  we deduce that $||df(x)||_{x}\le D_x^+f(x)$.

The proof of the second part follows along the same lines as  the Riemannian case (see \cite{GutuJaramillo}).
\end{proof}

\begin{rem}
It is worth noting that the above proof works for infinite-dimensional Finsler manifolds in the sense of Palais.
\end{rem}

In the next result we compare the upper scalar derivative with the norm of Clarke generalized differential for a locally Lipschitz mapping between Finsler manifolds. This relation will be of fundamental importance in order to obtain our global inversion results in Section 4.

\begin{lem}\label{Derjaco}
Let $M$ and $N$ be connected $C^1$ Finsler manifolds and let $f:M\to N$ a Lipschitz mapping on a neighborhood of $x\in M$. Then
$$
D_{x}^{+}f \leq |||\partial f(x)|||_{x}.
$$
\end{lem}
\begin{proof}
Let us fix $\varepsilon > 0$. Using  Lemma \ref{desigualdades:BiLipschitz} again, we can find a chart $(U,\varphi)$ of $M$ with  $x\in U$ such that $\varphi (U)$ is a convex set,  and a chart $(V,\psi)$ of $N$ with  $f(U)\subset V$, such that $\varphi$ and $\psi$ are $(1+\varepsilon)$-bi-Lipschitz for the norms $||\cdot||_{x,\varphi}$ on
$\mathbb R^m$ and $||\cdot||_{f(x),\psi}$ on $\mathbb R^n$, respectively. Then for every $y\in U$ we have that
\begin{align*}
d_{N}(f(x),f(y))  \leq & (1+\varepsilon)||\psi \circ f(x) - \psi \circ f(y)||_{f(x),\psi}\\
= & (1+\varepsilon) ||\psi\circ f \circ \varphi^{-1} (\varphi(x)) - \psi\circ f \circ \varphi^{-1} (\varphi(y))||_{f(x),\psi}.
\end{align*}
Applying the mean value theorem given in \cite[Proposition 2.6.5]{Clarkebook} to the mapping $\psi\circ f \circ \varphi^{-1}$, we obtain that
$$
\psi\circ f \circ \varphi^{-1} (\varphi(x)) - \psi\circ f \circ \varphi^{-1} (\varphi(y)) \in
\bigcup\limits_{w \in [\varphi(x),\varphi(y)]}\partial(\psi \circ f \circ \varphi^{-1})(w)(\varphi (x)-\varphi (y)).
$$
Then we have that

\begin{align*}
&d_{N}(f(x),f(y)) \leq  \\
& \leq (1+\varepsilon) \sup \left\{\| A \|_{x,\varphi,\psi}:  A \in \co\left(\bigcup\limits_{w \in [\varphi(x),\varphi(y)]}\partial(\psi \circ f \circ \varphi^{-1})(w)\right)\right\}||\varphi(x)-\varphi(y) ||_{x,\varphi}\\
&\le (1+\varepsilon)^{2}\sup\left\{\|A \|_{x,\varphi,\psi}:  A \in \co\left(\bigcup\limits_{w \in [\varphi(x),\varphi(y)]}\partial(\psi \circ f \circ \varphi^{-1})(w)\right)\right\}d_M(x,y).
\end{align*}

\noindent Since the mapping $\partial (\psi \circ f \circ \varphi^{-1})(\cdot)$ is upper semicontinuous, there exists
$\delta > 0$ such that, for all $w \in \overline{B}_{x,\varphi}(\varphi(x), \delta):=\{z\in \Real^m:
||z-\varphi(x)||_{x,\varphi} \leq \delta\}$, we get
\begin{equation*}
\partial (\psi\circ f \circ \varphi^{-1})(w) \subset \Jaco + \varepsilon \cdot \overline{B}_{x,\varphi,\psi},
\end{equation*}
where
$$
\overline{B}_{x,\varphi,\psi}:=\{A\in \mathcal{L}(\Real^m,\Real^n): ||A||_{x,\varphi,\psi} \leq 1\}.
$$
Therefore, using Proposition \ref{co-norm:igualdad},
\begin{equation*}
\| A \|_{x,\varphi,\psi}\leq \|\Jaco\|_{x,\varphi,\psi} + \varepsilon = |||\partial f(x)|||_x + \varepsilon,
\end{equation*}
for all $w\in \overline{B}_{x,\varphi}(\varphi(x),\delta)$ and
$A \in \partial (\psi\circ f \circ \varphi^{-1})(w)$.

\medskip

\noindent Now choose $\delta'>0$ such that $B_M(x,\delta^{\prime})\subset U$ and $\varphi(B_M(x,\delta')) \subset B_{x,\varphi}(\varphi(x),\delta)$. For every $y \in B_M(x,\delta^{\prime})$, we deduce that the set
$$
\bigcup\limits_{w \in [\varphi(x),\varphi(y)]}\partial(\psi \circ f \circ \varphi^{-1})(w)
$$
is contained into the ball $(|||\partial f(x)|||_x + \varepsilon)\cdot \overline{B}_{x,\varphi,\psi}\subset \mathcal{L}(\Real^m,\Real^n)$,  and therefore its convex hull is also contained into the same ball.

\

\noindent Thus, for every $y \in B_M(x,\delta^{\prime})$ we have that
$$
d_N(f(x),f(y)) \leq (1 + \varepsilon)^{2}(|||\partial f(x)|||_x + \varepsilon)d_M(x,y).
$$
As a consequence,
\begin{equation*}
D_{x}^{+}f = \lim\limits_{r \to 0}\left\{\sup\limits_{y\in B_M(x,r)}\frac{d_N(f(x),f(y))}{d_M(x,y)}\right\} \leq  (1+ \varepsilon)^{2}(|||\partial f(x)|||_{x} + \varepsilon),
\end{equation*}
for all $\varepsilon>0$, so that
$$
D_{x}^{+}f \leq |||\partial f(x)|||_x.
$$
\end{proof}


In the last part of this section, we are going to compare the norm of the Clarke generalized differential of a locally Lipschitz mapping $f$ between Finsler manifolds with the Lipschitz constant of $f$. What will provide a kind of mean value inequality in this context.


\begin{prop}\label{cotaja}
Let $M$ and $N$ be connected $C^1$ Finsler manifolds and  $f:M\to N$ a K-Lipschitz mapping on a neighborhood of $x\in M$. Then $|||\partial f(x)|||_x \leq K$.
\end{prop}

\begin{proof}
Fix  $\varepsilon > 0$. Using  Lemma \ref{desigualdades:BiLipschitz} we can find  charts $(U, \varphi)$ of $M$ at
$x$  and $(V, \psi)$  of $N$ at $f(x)$,  such that $f(U)\subset V$ and they are $(1+\varepsilon)$-bi-Lipschitz for the norms $||\cdot||_{x,\varphi}$ on $\mathbb R^m$ and $||\cdot||_{f(x),\psi}$ on $\mathbb R^n$, respectively. Thus the mapping
$\psi \circ f \circ \varphi^{-1}$ is $(1+\varepsilon)^2 K$-Lipschitz on a neighborhood of $\varphi(x)$ with the norms
$\|\cdot\|_{x,\varphi}$ and $\|\cdot\|_{f(x),\psi}$. Then it is not difficult to see that $\| A\|_{x,\varphi,\psi} \leq (1+\varepsilon)^{2} K$, for every $A\in \partial (\psi \circ f \circ \varphi^{-1})(\varphi (x))$.
Using  Proposition \ref{co-norm:igualdad} we obtain that
\begin{equation*}
|||\partial f(x)|||_{x} = \|\partial (\psi \circ f \circ \varphi^{-1})(\varphi (x))\|_{x,\varphi,\psi} \leq (1+\varepsilon)^{2} K.
\end{equation*}
Since the above inequality holds for every $\varepsilon>0$, we conclude that $|||\partial f(x)|||_x \leq K$.
\end{proof}

If $M$ is a connected Finsler manifold, the $\emph{length}$ of a continuous path $\gamma: [a,b] \to M$ is defined, as usual, as
\begin{equation*}
\ell(\gamma) = \sup\left\{\sum\limits_{i=1}^{k} d_M(\gamma(t_i), \gamma(t_{i-1})) \right\} \in [0, +\infty]
\end{equation*}
where the supremum is taken over all partitions $a=t_0< t_1<\cdots<t_k=b$. The path $\gamma:[a,b]\to M$ is  said to be \emph{rectifiable}  whenever $\ell(\gamma)<\infty$. In particular, every piecewise $C^1$-smooth
path $\gamma:[a,b]\rightarrow M$ is rectifiable, and in this case, as we remarked before, its length is given by the formula \eqref{lenght}.

\begin{prop}\label{rectilength}
Let $M$ and $N$ be connected $C^1$ Finsler manifolds and $f: M \to N$ a  locally Lipschitz mapping. If
$\gamma: [a, b] \to M$ is a rectifiable path and  $\sigma = f \circ \gamma:[a,b]\to N$, then
\begin{equation*}
\ell(\sigma) \leq \sup\left\{||| \partial f(z)|||_z \, : \,  z \in \Ima (\gamma)\right\}\ell(\gamma).
\end{equation*}
\end{prop}

\begin{proof}
The mean value inequality from  \cite[Proposition 3.8]{GutuJaramillo} gives us that
\begin{equation*}
\ell(\sigma) \leq \sup\left\{D^+ f(z) \, :\,  z \in \Ima (\gamma)\right\}\ell(\gamma),
\end{equation*}
so the result  follows from Lemma \ref{Derjaco}.
\end{proof}

\begin{cor}   Let $M$ and $N$ be connected $C^1$ Finsler manifolds
and let $f:M\to N$ be  a  locally Lipschitz  mapping. Then,  $f$ is Lipschitz if and only if  $||\partial f||_\infty := \sup\{|||\partial f(x)|||_x \, : \, x\in M\}<\infty$. Furthermore, $\Lip(f)= || \partial f||_\infty$.
\end{cor}
\begin{proof}
If $f$ is Lipschitz, Proposition \ref{cotaja} yields that $||\partial f||_\infty \leq \Lip(f)$. Conversely, suppose that $||\partial f||_\infty \leq K$. For each $x, y\in M$ and each $\varepsilon>0$, we can find a $C^1$-smooth
path $\gamma:[a,b]\rightarrow M$ from $x$ to $y$ such that $\ell (\gamma) \leq d_M(x, y)+ \varepsilon/K$. Then by Proposition \ref{rectilength} we have that
$$
d_N(f(x), f(y)) \leq \ell (f \circ \gamma) \leq K \, \ell (\gamma) \leq K d_M(x, y)+ \varepsilon.
$$
Hence, $f$ is $K$-Lipschitz.
\end{proof}
\section{Covering maps and global inversion}


We start this section with the problem of local inversion. For the case of locally Lipschitz mappings in $\mathbb R^n$,  Clarke obtained in \cite{Clarke1976} a local inversion theorem, and  Pourciau established in \cite{Pourciau1988}  the relation between the co-norm of the Clarke generalized differential of the map and the norm of the Clarke generalized differential of the corresponding local inverse. Now we are going to extend these results to the case of mappings between Finsler manifolds.

\begin{defn}
Let $M$ and $N$ be $n$-dimensional $C^1$ Finsler manifolds, and let $f:M\to N$ be a Lipschitz mapping on a neighborhood of $x\in M$.  We say that $\partial f(x)$ has \textbf{maximal rank} if  every element $B\in \partial f(x)$ has rank $n$, i.e., $B$ is an invertible linear map from $T_xM$ onto $T_{f(x)}N$.
\end{defn}

In what follows, if $\mathcal B $ is a family of invertible linear mappings, we denote by $\mathcal B^{-1}$
the set of their inverses $\mathcal B^{-1}=\{B^{-1} \, : \, B \in {\mathcal B} \}$, and we denote by $\co (\mathcal B^{-1})$ its convex hull.

\begin{prop}\label{thmfuninv}
Let $M$ and $N$ be $C^1$ Finsler manifolds of the same dimension, and let $f: M \to N$ be a Lipschitz mapping on a neighborhood of $x$  such that $\partial f(x)$ has maximal rank. Then there exist open neighborhoods $U$ of $x$ and $V$ of $f(x)$ such that $f_{\mid_U}: U \to V$ is a homeomorphism. Furthermore, if $g=(f_{\mid_U})^{-1}$, then $g$ is Lipschitz on a neighborhood of $f(x)$,   and we have that
$$
\partial g(f(x))\subset  \co (\partial f(x)^{-1})
$$
and
$$
|||\partial g(f(x))|||_{f(x)} \leq \frac{1}{\conormjax}.
$$
\end{prop}

\begin{proof}
Let $(U_0,\varphi)$ be a chart of $M$ at $x$ and $(V_0,\psi)$  a chart of $N$ at $f(x)$ such that $f$ is Lipschitz on $U_0$ and $f(U_0)\subset V_0$. Since $\psi \circ f \circ \varphi^{-1}$ is Lipschitz on a neighborhood of $\varphi(x)$ and $\Jaco$ has maximal rank, using the local inverse function theorem given by   Clarke in \cite{Clarke1976}, we obtain an open neighborhood $W$ of $\varphi(x)$ with  $W \subset \varphi(U_0)$ such that the restriction $(\psi \circ f \circ \varphi^{-1})_{\mid{W}}$ is a bi-Lipschitz homeomorphism onto its image. Setting $U=\varphi^{-1}(W)$ and $V=f(U)$, we have that $f_{\mid_U}: U\to V$ is a homeomorphism, and $g=(f_{\mid_U})^{-1}$ is Lipschitz on a neighborhood of $f(x)$.

\smallskip

\noindent On the other hand, since  $\varphi \circ g \circ \psi^{-1} = ( \psi \circ f_{\mid_U} \circ \varphi^{-1})^{-1}$, it follows from the results of  Pourciau in \cite[page 175]{Pourciau1988}  that
$$
\partial (\varphi \circ g \circ \psi^{-1}) (\psi (f(x)))\subset \co \left(\partial (\psi \circ f \circ \varphi^{-1})(\varphi (x))^{-1}\right).
$$
Then, it is clear that
$$
\partial g(f(x))\subset  \co (\partial f(x)^{-1}).
$$

\medskip

\noindent Let us prove the last inequality. Firstly,  notice that $\conormjax>0$ since every element in $\partial f(x)$ is invertible and $\partial f(x)$ is a compact subset
of ${\mathcal L}(T_xM, T_{f(x)}N)$. Now choose $B \in (\partial f(x))^{-1}$. Then $B^{-1} \in \partial f(x)$ and, by  \cite{Pourciau1988}  and Lemma \ref{norm:and:conorm},  it is easy to see that
\begin{equation}\label{igualnormas}
|||B|||_{f(x)}=  \frac{1}{/\hspace{-0.1cm}/\hspace{-0.1cm}/ B^{-1} /\hspace{-0.1cm}/\hspace{-0.1cm}/_{x}} \leq \frac{1}{\conormjax}.
\end{equation}
By convexity, the last inequality
also holds for every $B\in \co (\partial f(x)^{-1})$. Therefore
$$
|||\partial g(f(x))|||_{f(x)} \leq \frac{1}{\conormjax}.
$$
\end{proof}

Next we give our main results about global inversion. First recall that a continuous map between topological spaces $f:M\to N$ is said to be a \textbf{covering map} if every $z\in N$ has an open neighborhood $W$ such that $f^{-1}(W)$ is the disjoint union of open subsets of $M$ each of which is mapped homeomorphically onto $W$ by $f$.
Now, we can  use the co-norm of Clarke generalized differential in order to obtain an extension of Hadamard integral condition for locally Lipschitz mappings between Finsler manifolds.

\begin{thm}\label{main:thm}
Let $M$ and $N$ be connected $C^{1}$ Finsler manifolds of the same dimension, where $M$ is complete, and
let  $f: M \to N$ be a locally Lipschitz mapping such that $\partial f(x)$  has maximal rank for all $x\in M$.
Assume that there exists $x_0 \in M$ such that
\begin{equation*}
\int_{0}^{\infty} m(t)dt = \infty, \quad \text{where} \quad m(t) = \inf\limits_{x \in \overline{B}_M(x_0,t)}/\hspace{-0.1cm}/\hspace{-0.1cm}/ \partial f(x)/\hspace{-0.1cm}/\hspace{-0.1cm}/_{x}.
\end{equation*}
Then $f$ is a covering map.
\end{thm}

\begin{proof}
First of all, notice that $m(t)>0$ for all $t>0$. Indeed, if  there is $t_0>0$ such that $m(t_0)=0$, then $0\le m(t)\le m(t_0) =0$ for all $t\ge t_0$ and, since $m(t)\le /\hspace{-0.1cm}/\hspace{-0.1cm}/ \partial f(x_0)/\hspace{-0.1cm}/\hspace{-0.1cm}/_{x_0}$ for all $t>0$, we obtain that
\begin{equation*}
\int_{0}^{\infty} m(t)dt = \int_0^{t_0} m(t)dt \le t_0  /\hspace{-0.1cm}/\hspace{-0.1cm}/ \partial f(x_0)/\hspace{-0.1cm}/\hspace{-0.1cm}/_{x_0} < \infty,
\end{equation*}
which contradicts the hypothesis.

 On the other hand, for every $x \in M$, we know from Proposition \ref{thmfuninv} that  $f$ is a local homeomorphism around $x$ and that, if $g$ is the local inverse of $f$ at $f(x)$, using Lemma \ref{Derjaco} we have that
$$
D_{f(x)}^{+} g \leq |||\partial g(f(x))|||_{f(x)} \leq \frac{1}{\conormjax}.
$$
Thus
$$
D_{x}^- f = \frac{1}{D_{f(x)}^{+} g } \geq \conormjax.
$$
\noindent Now let us define $w(t) = \frac{1}{m(t)}$. It is clear that $w$ is a weight, i.e., $w: (0,\infty) \to [0,\infty)$ is a nondecreasing map (not necessarily continuous) such that $\displaystyle{\int\limits_{0}^{\infty} \frac{dt}{w(t)} = \infty}$.

\noindent Since $f$ is a local homeomorphism, and using   \cite[Lemma 4.5 and Theorem 4.6]{GutuJaramillo}, it  is sufficient to prove that $D_{x}^{-}f\cdot w(d_M(x,x_0)) \geq 1$, for every $x \in M$. And taking into account that
$x\in \overline{B}_M(x_0, d_M(x, x_0))$,   we can finish the proof:
\begin{align*}
D_{x}^{-} f \cdot w(d_M(x,x_0)) \geq & /\hspace{-0.1cm}/\hspace{-0.1cm}/ \partial f(x)/\hspace{-0.1cm}/\hspace{-0.1cm}/_{x} \, \frac{1}{m(d_M(x,x_0))} \\
  =  & \frac{/\hspace{-0.1cm}/\hspace{-0.1cm}/ \partial f(x)/\hspace{-0.1cm}/\hspace{-0.1cm}/_{x}}{\inf\limits_{z \in \overline{B}_M(x_0, d_M(x, x_0))} /\hspace{-0.1cm}/\hspace{-0.1cm}/ \partial f(z)/\hspace{-0.1cm}/\hspace{-0.1cm}/_{z}} \ge \frac{/\hspace{-0.1cm}/\hspace{-0.1cm}/ \partial f(x)/\hspace{-0.1cm}/\hspace{-0.1cm}/_{x}}{/\hspace{-0.1cm}/\hspace{-0.1cm}/ \partial f(x)/\hspace{-0.1cm}/\hspace{-0.1cm}/_{x}} =  1.
\end{align*}

\end{proof}

Let $f:M\to N$ be a covering map between path-connected metric spaces, and $f_*:\pi_1(M)\to \pi_1(N)$  the associated morphism between their fundamental groups. It is well known that $f$ is a homeomorphism from $M$ onto $N$ if, and only if, $f_*[\pi_1(M)]=\pi_1(N)$ (see e.g. \cite[Chapter 2]{Spanier}). Thus we obtain the following corollary.

\begin{cor}
Under the assumptions of Theorem \ref{main:thm}, assume that either $N$ is simply connected or $\pi_1(M)=\pi_1(N)$ is finite. Then $f$ is a global homeomorphism.
\end{cor}


In the last part of the paper we use a variant of Hadamard integral condition to obtain global inversion and global injectivity results for Lipschitz mappings
in $\mathbb R^n$ in terms of spectral conditions of the Clarke generalized differential. In particular, we will give  extensions of some of the results given in  \cite{BiGuSan} and \cite{FeGuRa} to a nonsmooth setting.
For a locally Lipschitz mapping  $f: \Real^n \to \Real^n$, we denote by $\Spec(f)$ the set
of all complex eigenvalues of all matrices $A \in \partial f(x)$ for all $x \in \Real^n$.
\begin{thm}
Let $f: \Real^n \to \Real^n$ be a Lipschitz mapping such that $\partial f(x)$ has maximal rank for every $x \in \Real^n$. Assume that there exists $x_0 \in \Real^n$ such that
\begin{equation*}
\int_{0}^{\infty} s(t)dt = \infty, \quad \text{where} \quad s(t) = \inf\limits_{x \in \overline{B}(x_{0},t)} \left\{ |\lambda|^{n}: \lambda \in \Spec(\partial f(x)) \right\}.
\end{equation*}
Then $f$ is a global homeomorphism.
\end{thm}

\begin{proof}
Suppose that $f$ is $K$-Lipschitz. By Proposition \ref{cotaja} we have that $\|A\|  \leq K$ for every $A\in \partial f(x)$ and every $x \in \Real^n$. Let us denote by $\|\cdot\|_\infty$ the supremum norm on the space of real $n\times n$ matrices, i.e.
\begin{equation*}
\|A\|_\infty:= \sup\{|a_{i,j}|:1\leq i,j \leq n\} \quad  \text{ where } \quad  A= (a_{ij}).
\end{equation*}
Since $\|\cdot\|$  and $\|\cdot\|_{\infty}$ are equivalent norms, there exists a  constant $K_{1}>0$ such that
\begin{equation*}
\|A\|_{\infty} \leq K_{1} \quad \text{for every $A\in \partial f(x)$ and every $x\in \Real^n$.}
\end{equation*}
Therefore,  $\|\adj (A)^t\|_{\infty} \leq (n-1)! {K_{1}}^{n-1} = K_{2}$ for every $A\in \partial f(x)$ and $x\in \Real^n$, where $\adj (A)^t$ denotes the  transpose of the    adjoint matrix of $A$.

\smallskip

\noindent On the other hand, fixed  $t>0$, $x\in \overline{B}(x_0,t)$  and $A\in \partial f(x)$,
we have that  $\det A = \prod\limits_{i = 1}^{n} \lambda_i$ where $\lambda_i$ are the eigenvalues of $A$, so that
\begin{align*}
|\det A|  & =  |\lambda_1| \cdotÊ |\lambda_2|\cdotsÊ|\lambda_n| \ge \inf\{|\lambda|^{n}: \lambda \in \Spec(A)\}   \\
&  \ge  \inf\limits_{x \in \overline{B}(x_{0},t)}\{|\lambda|^{n}: \lambda \in \Spec(\partial f(x))\} =s(t).
\end{align*}

\noindent Hence, we have that the inverse matrix $A^{-1} = \frac{1}{\det A}\, \adj (A)^t$ satisfies
\begin{equation*}
\|A^{-1}\|_{\infty} \leq \frac{1}{|\det A|}\|\adj (A)^t\|_{\infty} \leq \frac{1}{|\det A|} K_2 \leq  \frac{K_2}{s(t)}.
\end{equation*}
Using again that  the norms $||\cdot||$ and $||\cdot||_\infty$ are equivalent, we get a constant $K_3>0$ such that
\begin{equation*}
\|A^{-1}\|^{-1} \ge K_3 s(t) \quad \text{for every $A\in \partial f(x)$ and $x \in \overline{B}(x_0,t)$.}
\end{equation*}
Therefore
\begin{align*}
K_3 s(t) \le &  \inf\limits_{x\in\overline{B}(x_0,t)} \left\{ \inf\limits_{A \in \partial f(x)} \|A^{-1}\|^{-1} \right\} \\
& =  \inf\limits_{x\in\overline{B}(x_0,t)} \left\{ \inf\limits_{A \in \partial f(x)}  /\hspace{-0.1cm}/ A /\hspace{-0.1cm}/ \right\} \\
& = \inf\limits_{x\in\overline{B}(x_0,t)}/\hspace{-0.1cm}/\partial f(x) /\hspace{-0.1cm}/ = m(t).
\end{align*}
Hence, $\displaystyle{\int\limits_{0}^{\infty}m(t)dt \geq \int\limits_{0}^{\infty}K_3  s(t)dt = \infty}$ and by Theorem \ref{main:thm}  (see also \cite{Pourciau1988}), the mapping $f$ is a global homeomorphism.
\end{proof}

As a direct consequence of the previous theorem, we obtain the following corollary.

\begin{cor}\label{cor:spec}
Let $f: \Real^n \to \Real^n$ be a Lipschitz mapping such that $\partial f(x)$ has maximal rank for all $x \in \Real^n$,
and suppose that there is some $\varepsilon>0$ satisfying
\begin{equation*}
\Spec(\partial f(x)) \cap \{z \in \C: |z| \leq \varepsilon\} = \emptyset, \quad \text{ for all $x\in \Real^n$.}
\end{equation*}
 Then $f$ is a global homeomorphism.
\end{cor}

Finally, we are going to give a global injectivity result for nonsmooth mappings. First we need the following lemma, which is a nonsmooth version of Main Lemma in \cite{FeGuRa}.

\begin{lem}\label{injective}
Let $f: \Real^n \to \Real^n$ be a Lipschitz mapping such that $\partial f(x)$ has maximal rank for every $x \in \Real^n$. Given $t\in \mathbb R$, consider the map $f_t(x):=f(x) - tx$. If there exists a sequence $(t_k)$ of real number converging to $0$ such that every map $f_{t_k}$ is injective, then $f$ is injective.
\end{lem}

\begin{proof}
Let $x_1, x_2 \in \mathbb R^n$ such that $f(x_1)=f(x_2)=y$. We are going to show that $x_1=x_2$. Let us define the mapping $h: \mathbb R \times \mathbb R^{n} \to \mathbb R \times \mathbb R^{n}$ by $h(t, x)= (t, f(x)-tx)$. Then $h$ is locally Lipschitz, and it is not difficult to see that, for every $x\in \mathbb R^n$, $h$ has maximal rank at $(0, x)$.  Applying the local inverse function theorem from \cite{Clarke1976} to $h$ at $(0, x_1)$ and $(0, x_2)$, we obtain open neighborhoods $W_1$ of $(0, x_1)$,  $W_2$ of $(0, x_2)$ and  $V$ of $y$, and some $\delta>0$ such that $h_{\mid_{W_1}}: W_1 \to (-\delta, \delta) \times V$ and $h_{\mid_{W_2}}: W_2 \to (-\delta, \delta) \times V$ are homeomorphisms. For $i=1, 2$ set $g_i= (h_{\mid_{W_i}})^{-1}$,  and note that $g_i(0, y)=(0, x_i)$ for $i=1,2$. Now $t_k \in (-\delta, \delta)$ for $k$ large enough, and for $i=1, 2$ we have that $g_i(t_k, y)=(t_k, x_k^i)$ for some $x_k^{i}\in W_i$ satisfying $y=f_{t_k}(x_k^{i})$. Since each $f_{t_k}$ is injective, we deduce that $x_k^{1}=x_k^{2}$ for $k$ large enough. But for $i=1, 2$, we know that $x_i= \lim_{k \to \infty}x_k^{i}$. In this way we obtain that $x_1=x_2$.
\end{proof}


\begin{thm}
Let $f: \Real^n \to \Real^n$ be a Lipschitz mapping such that $\partial f(x)$ has maximal rank for every $x \in \Real^n$. Suppose that there exists a sequence $(D_k)$ of compact discs of $\mathbb C$ (with nonempty interior), centered at points $t_k$ of the real axis, such that $\lim_{k \to \infty} t_k =0$ and
\begin{equation*}
\Spec(\partial f(x)) \cap (\cup_{k=1}^{\infty} D_k) = \emptyset.
\end{equation*}
Then $f$ is injective
\end{thm}
\begin{proof}
Let us denote the mappings $f_m(x)= f(x)-t_mx$ for all $t_m$. It is clear that $\Spec(\partial f_m(x)) \cap (D_m-t_m) = \emptyset$, where $D_m -t_m = \{z \in  \C : z+t_m \in D_m\}$ is a compact disc centered at $0$. Using Corollary \ref{cor:spec}, we obtain that $f_m$ is a global homeomorphism for every $m$. Lemma \ref{injective} allows us to conclude that $f$ is injective.
\end{proof}

Finally, as a consequence, we obtain:

\begin{cor}
Let $f: \Real^n \to \Real^n$ be a Lipschitz mapping such that $\partial f(x)$ has maximal rank for all $x \in \Real^n$,
and suppose that
$$
\Spec(\partial f(x)) \subset \{z \in \C \, : \, {\mathfrak R}(z) <0 \}.
$$
Then $f$ is injective.
\end{cor}


\end{document}